\documentclass[final,3p,11pt, times]{article}
\usepackage{authblk}
\usepackage{amsmath}
\usepackage{amssymb,latexsym}
\usepackage[german,english]{babel}
\usepackage{url}
\usepackage{graphicx}
\usepackage{gastex}
\usepackage{longtable}
\usepackage{lscape}
\usepackage{tabularx}
\usepackage{multicol}
\usepackage{verbatim}
\usepackage{multirow}
\usepackage{graphicx}
\usepackage{epsfig,graphics,graphicx}
\usepackage{geometry}
\usepackage[usenames]{xcolor}

\usepackage{tikz}

\usepackage{booktabs}
\usepackage{hhline}

\hbadness=\maxdimen

\newtheorem{lem}{{\bf Lemma}}[subsection]
\newtheorem{thm}{{\bf Theorem}}[section]

\usepackage{chngcntr}
\counterwithout{equation}{section}

\begin{document}
	
	\title{Hosoya properties of power graphs over certain groups}
	
	\author[1]{Yogendra Singh}
	\author[1]{Anand Kumar Tiwari}
	\author[2]{Fawad Ali}
	
	\affil[1]{\small Department of Applied Science, Indian Institute of Information Technology, Allahabad 211015, India}

	\affil[2]{Institute of Numerical Sciences, Kohat University of Science and Technology, Kohat 26000, KPK, Pakistan \hspace{15cm} E-mail: yogendraiiita@gmail.com, anand@iiita.ac.in, fawadali@math.qau.edu.pk}

	\maketitle
	
	\hrule
	
\begin{abstract}
	
	The power graph denoted by $\mathcal{P}(\mathcal{G})$ of a finite group $\mathcal{G}$ is a graph with vertex set $\mathcal{G}$ and there is an edge between two distinct elements $u, v \in \mathcal{G}$ if and only if $u^m = v$ or $v^m = u$ for some $m \in \mathbb{N}$. Depending on the distance, the Hosoya polynomial contains a lot of knowledge about graph invariants which can be used to determine well-known chemical descriptors. The Hosoya index of a graph $\Gamma$ is the total number of matchings in $\Gamma$. In this article, the Hosoya properties of the power graphs associated with a finite group, including the Hosoya index, Hosoya polynomial, and its reciprocal are calculated.
	  
\end{abstract}


\textbf{Keywords:} Power graph, Hosoya index, Hosoya polynomial.

\textbf{MSC(2010):} 05C50, 05C07, 15A27.

\hrule

\section{Introduction}

A topological index is a number derived from the configuration of a molecule that indicates its fundamental structural characteristics.
A graph topological index may be classified into several types, the most notable of which are distance-based, degree-based, eccentricity-based, and spectrum-based topological indices. 

 As discussed above, topological indices are used to investigate several properties of chemical compound with the help of its molecular structure (graph). In this regards, Wiener \cite{Wie(1947)} introduced the first topological index which is known as path number. Further, various researchers utilized P\'{o}lya's \cite{polya1937kombinatorische} concept of evaluating polynomials to determine the unsaturated hydrocarbons molecular orbital. This idea was used by Hosoya \cite{hosoya1988some} in $1988$ to establish polynomials of various chemical structures that is referred to as the \emph{Hosoya polynomials} and got attention from all over the world.
This polynomial was called the  \emph{Wiener polynomial} by Sagan et al. \cite{sagan1996wiener} in $1996$, although according to the majority of academics, it is still known as the  \emph{Hosoya polynomial}.

In this article, we consider undirected simple graphs. 
Kelarev and Quinn \cite{kelarav(2002)} introduced the notion of directed power graphs over semigroups having vertex set $S$ and there is an edge from the vertex $x$ to $y$ iff $x \neq y$ and $y = x^m$ for some $m \in \mathbb{N}$. Motivated by this, Chakrabarty et al. \cite{1} constructed an undirected power graph for a group in which two distinct vertices are connected by an edge if and only if one can be expressed as an integral power of the other.
Recently, the power graph has become an important topic of research in several areas of mathematics, such as group theory, ring theory, and gyrogoups; see \cite{ys(2022), stap(2022), sta(2022)}. The power graph was first reviewed in $2013$ \cite{survey0}, and the most recent review was in $2021$ \cite{kumar2021recent}.

Different types of topological indices have been calculated for (non) commuting graphs in \cite{Wei(2020), arsu(2022),  Ali2022} and for power graphs over finite groups in \cite{rather(2022), ali(2022)}. Following their work, we focus on the power graphs of a finite group. In this paper, we evaluate various Hosaya polynomials and Hosaya index of the power graph of a finite group.

The remaining sections are arranged as follows: Section \ref{Sec2} provides definitions, symbols, and findings that will be required later. In Section \ref{Sec3}, Hosoya polynomials and reciprocal Hosoya polynomials of the power graph of a group is determined. In Section \ref{Sec4}, Hosoya index of the power graph of a group is determined.

\section{Preliminaries} \label{Sec2}

This section covers some fundamental graph-theory concepts and well-known results that will be useful throughout the article, such as \cite{cha(2006), West(2009)}.

Let $\Gamma$ be a simple and undirected graph. Let $V(\Gamma)$ and $E(\Gamma)$ represent the vertex and edge sets of $\Gamma$, respectively. If two vertices $v_1$ and $v_2$ are connected by an edge, we write as $v_1 \sim v_2$, otherwise as $v_1 \nsim v_2$. The degree of a vertex $v$, denoted by $\deg(v)$, is the number of vertices connected with $v$. The shortest distance between $u_1$ and $u_2$ in $\Gamma$ is the distance from $u_1$ to $u_2$, which is represented by $d(u_1$, $u_2$). The largest distance between a vertex $u$ to all the other vertices of $\Gamma$ is known as the eccentricity of $u$ and is denoted by the symbol $ec(u)$. The radius of $\Gamma$, denoted by $rad(\Gamma)$, is the smallest eccentricity among all the vertices of $\Gamma$. The diameter of $dia(\Gamma)$, denoted by $dia(\Gamma)$, is the largest eccentricity among all the vertices of $\Gamma$.
The join of two connected graphs $\Gamma_{1}$ and $\Gamma_{2}$, denoted by $\Gamma_{1} \vee \Gamma_{2}$, is a graph with the vertex and edge sets $ V(\Gamma_{1}) \cup V(\Gamma_{2})$ and $ E(\Gamma_{1}) \cup E(\Gamma_{2}) \cup \{ c \sim d : c \in V(\Gamma_{1}),~ d \in V(\Gamma_{2})\}$, respectively. Let $K_n$ denote complete graph on $n$ vertices. For a graph $\Gamma$, Hosoya defined the following polynomial, also named as Hosoya Polynomial.
\begin{equation}\label{definition of hosaya poly}
\mathcal{H}(\Gamma, x)= \sum_{i \geq 0} dis(\Gamma, i) x^{i}.
\end{equation}
The coefficient $dis(\Gamma, i)$ denotes the total number of $(v, w)$ pairs of vertices such that $dis(v, w)=i,$ where $i \leq diam(\Gamma)$.
In \cite{Ramane(2019)}, the authors presented the reciprocal status Hosoya polynomial for $\Gamma$ as

\begin{equation}\label{definition of reciprocal hosaya poly}
\mathcal{H}_{rs}(\Gamma, x)= \sum_{vw \in E(\Gamma)} x^{rs(v)+rs(w)},
\end{equation}
where $rs(w)= \sum\limits_{v \in V(\Gamma), v \neq w}\frac{1}{dis(w, v)}$ is called the reciprocal status or the transmission of $w$.

This article focuses on the group $G^{2^{k-1}p-1}_{2(2^kp)}$ with the following representation
 $$G^{2^{k-1}p-1}_{2(2^kp)} = \langle r, s \, : r^{2^kp} = s^2 = e,~ srs^{-1} = r^{2^{k-1}p-1}\rangle,$$ where $G^{2^{k-1}p-1}_{2(2^kp)}$ is a non-abelian group of order $2(2^kp)$, where  $k \geq 2$ and $p$ is an odd prime.
 
Throughout this paper, we denote the group $G^{2^{k-1}p-1}_{2(2^kp)}$ by $\mathcal{G}$. We consider the following partition of the group $\mathcal{G} = \{H_0 , H_1, H_2, H_3 \}$, where $H_0 = \{e, u = r^{2^{k-1}p}\},~
H_1 =  \{r^i : 1 \leq i \leq 2^kp-1\} \ \setminus \ \{u\},~
H_2 =  \{sr^{2t} {: 1 \leq t \leq 2^{k-1}p}\},$ and $H_3 = \{sr^{2j+1} : 0 \leq j \leq 2^{k-1}p-1\}$.

	For the above structure, we have

\begin{thm}
	For the group $\mathcal{G}$, we have
	$$P(\mathcal{G}) \cong P \left(\mathbb{Z}_{2^kp} \right) \cup 2^{k-1}p K_2 \cup 2^{k-2}p \left((P \langle r^{2^{k-1}p} \rangle) + K_2 \right).$$
	
\end{thm}

\section{Hosoya properties}\label{Sec3} 

In this section, we discuss the Hosoya polynomial and reciprocal status Hosoya polynomial for the power graph of group $\mathcal{G}$. 

\begin{lem} The Hosoya polynomial $H (P(\mathcal{G}), x )$ of the power graph $P(\mathcal{G})$ is 
	
	\[H (P(\mathcal{G}), x ) = \left(3.2^{2k-1}p^2-9.2^{k-2}p \right)x^2 + \left(2^{2k-1}p^2 + 5p.2^{k-2} \right) x + 2^{k+1}p.\]
\end{lem}

\noindent{\bf Proof.} 
Let $V_p$ be the set of all pairs of vertices of $P(\mathcal{G})$. 
Then 
$$|V_p| = 2^{k+1}p + \frac{2^{k+1}p(2^{k+1}p-1)}{2}.$$ Now, suppose 
$$C(P(\mathcal{G}),i)= \left\{(l,m); l, m \in V(P(\mathcal{G})) : dis(l,m) = i \right\}.$$ 
Then $ dis(P(\mathcal{G}),i) = |C(P(\mathcal{G}),i)|$. As $diam(P(\mathcal{G})) = 2$, we need to find $dis(P(\mathcal{G}), 0)$,
$dis(P(\mathcal{G}), 1)$, and $dis(P(\mathcal{G}), 2)$. As, $d(l,l) = 0  \text{ for all } l \in V(P(\mathcal{G}))$, we have 
$|C(P(\mathcal{G}), 0)| = 2^{k+1}p.$ 

Now, 
$C(P(\mathcal{G}),1) = \left\{(e,h_1): h_1 \in H_1 \cup \{u\} \right\} \cup \left\{(e,h_2): h_2 \in H_2\right\} \cup  \left\{(e,h_3):  h_3 \in H_3 \right\} \cup \left\{(u,h_3): h_3 \in H_3 \right\} \cup \left\{(y,z); y,z \in H_3 : y = rs^{2j+1}, z = rs^{2j+1+2^{k-1}p} \text{ for } 0 \leq j \leq 2^{k-2}p-1 \right\}$.
This implies 
\begin{align*}
|C(P(\mathcal{G}),1)| & = dis(P(\mathcal{G}), 1)\\
& = \frac{2^kp(2^kp-1)}{2} + 2^{k-1}p + 2^{k-1}p  + 2^{k-1}p + 2^{k-2}p \\ & =  2^{2k-1}p^2 + 5p.2^{k-2}.
\end{align*}
Clearly, $$V_p = C(P(\mathcal{G}), 0) \cup C(P(\mathcal{G}), 1) \cup C(P(\mathcal{G}), 2).$$
Thus, we have $$|V_p| = dis(P(\mathcal{G}), 0) + dis(P(\mathcal{G}), 1) + dis(P(\mathcal{G}), 2).$$

Hence, 
\begin{equation*} 
\begin{split}
dis(P(\mathcal{G}), 2) & = |V_p| - dis(P(\mathcal{G}), 0) - dis(P(\mathcal{G}), 1)  \\
& = 2^{k+1}p + \frac{2^{k+1}p(2^{k+1}p-1)}{2} - 2^{k+1}p - 2^{2k-1}p^2 - 5p.2^{k-2}\\
& = 3.2^{2k-1}p^2 - 9.2^{k-2}p.
\end{split}
\end{equation*} 

Since $$H(P(\mathcal{G}), x) = \sum_{i \geq 0} dis(\mathcal{G}, i)x^i, \ \text{where} \ i \leq diam(\mathcal{G})=2 \text{, we have}$$ 
$$H (P(\mathcal{G}), x ) = (3.2^{2k-1}p^2-9.2^{k-2}p)x^2 + (2^{2k-1}p^2 + 5p.2^{k-2}) x + 2^{k+1}p.$$ \hfill $\Box$ 

\begin{lem} The reciprocal Hosoya polynomial 	$H_{rs}(P(\mathcal{G})), x)$ of the graph $P(\mathcal{G})$ is
	\begin{align*}
	H_{rs}(P(\mathcal{G})), x) &= 
	x^{15.2^{k-2}p-2} + (2^kp-2).x^{7.2^{k-1}p-2} + 2^{k-1}p.x^{3.2^{k}p-2} + 2^{k-1}p.x^{3.2^kp} \\
	&+ 2^{k-1}p.x^{11.2^{k-2}p} + \frac{(2^{k}p-2)(2^{k}p-2-1)}{2}.x^{3.2^kp-2} +
	2^{k-2}p.x^{2^{k+1}p+2}.
	\end{align*}
\end{lem}	
\noindent{\bf Proof.} Observe that, there are seven types of edges in $P(\mathcal{G})$, namely $eu$, $eh_1$, $eh_2$, $eh_3$, $uh_3$, $vw$, $yz$, where $v,w \in H_1$, $h_i \in H_i$ for $1 \leq i \leq 3$, and $y=rs^{2j+1}$, $z = rs^{2j+1+2^{k-1}p} \in H_3$ for $ 0 \leq j \leq 2^{k-2}p-1$. Then by Equation (\ref{definition of reciprocal hosaya poly}), we have
\hspace{-2cm}
\begin{align*}
H_{rs}(P(\mathcal{G}), x) & =    \sum\limits_{eh_1 \in E(G)} x^{rs(e)+rs(h_1)} + \sum\limits_{eh_2 \in E(G)} x^{rs(e)+rs(h_2)}  + \sum\limits_{eh_3 \in E(G)} x^{rs(e)+rs(h_3)} \ + \\ &  \sum\limits_{uh_3 \in E(G)} x^{rs(u)+rs(h_3)} + \sum\limits_{vw \in E(G)}x^{rs(v)+rs(w)} + \sum\limits_{yz \in E(G)}x^{rs(y)+rs(z)} + x^{rs(e)+rs(u)}.
\end{align*}

Note that,
$deg(e)=2^{k+1}p-1$, $deg(u)=3.2^{k-1}p-1$, $deg(h_1)=2^kp-1$, $deg(h_2)=1$, $deg(h_3)=3$. 
As $diam(P(\mathcal{G}))=2$, and for any vertex $s \in P(\mathcal{G})$, we have $deg(s)$ number of vertices at distance $1$ form $s$ and $2^{k+1}p-1-deg(s)$ vertices at distance $2$. 
Hence, $rs(s) =\frac{1}{1}.deg(s) + \frac{1}{2}.(2^{k+1}p-1-deg(s))$. Using this, we have
\begin{align*}
rs(e) &=  \frac{1}{1}.(2^{k+1}p-1) + \frac{1}{2}.(2^{k+1}p-1-2^{k+1}p+1) = 2^{k+1}p-1,\\
rs(h_1) &= rs(v) = rs(w) = \frac{1}{1}.(2^{k}p-1) + \frac{1}{2}.(2^{k+1}p-1-2^{k}p+1) = 2^{k+1}p-2^{k-1}p-1,\\
rs(h_2) &= \frac{1}{1}.(1) + \frac{1}{2}.(2^{k+1}p-1-1) = 2^kp,\\
rs(h_3) &= rs(y) = rs(z) = \frac{1}{1}.(3) + \frac{1}{2}.(2^{k+1}p-1-3) = 2^{k}p+1,\\
rs(u) &= \frac{1}{1}.(3.2^{k-1}p-1) + \frac{1}{2}.(2^{k+1}p-1-3.2^{k-1}p+1) = 7.2^{k-2}p-1.
\end{align*}
Thus, $rs(e)+rs(h_1) = 7.2^{k-1}p-2, \ \  rs(e)+rs(h_2) = 3.2^{k}p-1, \ \ rs(e)+rs(h_3) = 3.2^{k}p,$ \\
$rs(u) + rs(h_3)  = 11.2^{k-2}p, \ \ rs(v)+rs(w) = 3.2^{k}p-2,  \ \ rs(y)+rs(z) = 2^{k+1}p+2, rs(e)+rs(u) = 15.2^{k-2}p-2.$ 

\smallskip
Therefore, the number of edges in $P(\mathcal{G})$ of types $eu$, $eh_1$, $eh_2$, $eh_3$, $uh_3$, $vw$, $yz$ are 1, $2^{k}p-2, 2^{k-1}p, 2^{k-1}p, 2^{k-1}p,  \frac{(2^{k}p-2)(2^{k}p-2-1)}{2}, 2^{k-2}p$, respectively. By using these values in $H_{rs}(P(\mathcal{G}), x)$, 
we get the required proof. \hfill $\Box$

\section{Hosaya Index} \label{Sec4} 

This section examines the Hosoya index of the power graph of group $\mathcal{G}$. The number of non-void matchings of order $m_i$ in $K_{2^kp}$ is given in Table 1, where $1 \leq i \leq \frac{2^kp}{2}.$

\begin{center}
	
	\noindent \textbf{Table 1} : The total number of non-void matchings in $K_{2^kp}$
	
	\renewcommand{\arraystretch}{2.5}
	\begin{tabular}{|p{1.9cm}|p{.75cm}|p{2.2cm}|p{3.4cm}|p{.8cm}|p{3.2cm}|}
		
		\hline
		Order & $m_{1}$ & $m_{2}$ & $m_{3}$  &~~~ $\cdots$ & $m_{i} \ (1 \leq i \leq 2^{k-1}p)$ \\
		\hline
	Number of matchings	& $\binom{2^kp}{2}$ & $\frac{1}{2}\binom{2^kp}{2}\binom{2^kp-2}{2}$ & $\frac{1}{3}\binom{2^kp}{2}\binom{2^kp-2}{2}\binom{2^kp-4}{2}$  &~~~ $\cdots$ & $\frac{1}{i} \prod_{s=0}^{i-1}\binom{2^kp-2s}{2}$ \\
		
		\hline
	\end{tabular}
\end{center}

\smallskip
\begin{thm} \label{theorem1}
	
	The Hosoya index of $\mathcal{P}(\mathcal{G})$ is given by 
	
	$1 + \sum\limits_{i=1}^{2^{k-1}p} M_{1}^{i} \ + \ M^1_2 + \sum\limits_{i=1}^{2} M_{3}^{i} + 
	\sum\limits_{i=1}^{2^{k-2}p} M_{4}^{i} + \sum\limits_{i=2}^{2^{k-1}p+1} M_{5}^{i} \ + \
	\sum\limits_{i=2}^{3 \cdot 2^{k-2}p} M_{6}^{i} + M^{2}_7 + \sum\limits_{i=3}^{2^{k-2}p-1} M_{7}^{i} + 
	\sum\limits_{i=2}^{2^{k-1}p} M_{8}^{i} + M^2_9 + \sum\limits_{i=2}^{2^{k-2}p+1} M_{10}^{i} + \sum\limits_{i=3}^{3 \cdot 2^{k-2}p} M_{11}^{i} +
	\sum\limits_{i=3}^{3 \cdot 2^{k-2}p} M_{12}^{i} + \sum\limits_{i=3}^{2^{k-2}p+1} M_{13}^{i}  + \sum\limits_{i=3}^{ 2^{k-1}p+1} M_{14}^{i} +  + \sum\limits_{i=4}^{3 \cdot 2^{k-2}p} M_{15}^{i}  ,$ where for $ 1 \leq j \leq 15$ $M^{i}_{j}$ is defined explicitly in the proof.

\end{thm}

\noindent{\textbf{Proof.}} Consider the subsets 
$ A_1 = \{e,r,r^2, \ldots, r^{2^kp-1}\}, \
A_2 = \{sr^{2t} {: 1 \leq t \leq 2^{k-1}p}\}, \ A_3 = A_1 \setminus \{e\}, \
A_4 = A_1 \setminus \{r^{2^{k-1}p}\}, \ \Omega = \{e,r^{2^{k-1}p}\}, \
A_5 = A_1 \setminus \Omega, \
A_6 = \bigcup_{j=0}^{2^{k-2}p-1} A^j_6, \text{ where } A^j_6 = \{sr^{2j+1}, sr^{2j+1+2^{k-1}p}\}$ of $\mathcal{G}$. Using this, we have the following types of edges in $\mathcal{P}(\mathcal{G}).$

\begin{center}
	
	\noindent \textbf{Table 2} : Types of edges
	
	\renewcommand{\arraystretch}{.9}
	\begin{tabular}{|p{6cm}|p{7.5cm}|p{8cm}|}
		
		\hline
		
		\hline
		Type &  Type   \vspace{.25cm} \\  \hline
		 E-$1:$ $x_1 \sim x_2$, for $x_1, x_2 \in A_{1}$ & E-$8:$ $x_1 \sim x_2$, for $x_1, x_2 \in \Omega$ \vspace{.25cm} \\ 
	     E-$2:$ $x_1 \sim x_2$, for $x_1 \in A_{1}, x_2 \in A_6$ & E-$9:$ $x_1 \sim x_2$, for $x_1 \in A_{5},  x_2 \in \Omega$ \vspace{.25cm} \\ 
		 E-$3:$ $x_1 \sim x_2$, for $x_1 \in A_{1}, x_2 \in A_2$ & E-$10:$ $x_1 \sim x_2$ for $x_1 \in A_{6},  x_2 \in \Omega$ \vspace{.25cm} \\ 
	     E-$4:$ $x_1 \sim x_2$, for $x_1, x_2 \in A_{3}$ & E-$11:$ $x_1 \sim x_2$, for $x_{1},  x_{2} \in A_{6}^{j} \subseteq A_{6}$ \vspace{.25cm} \\
		 E-$5:$ $x_1 \sim x_2$, for $x_1 \in A_{3}, x_2=e$ & E-$12:$ $x_1 \sim x_2$, for $x_{1},  x_{2} \in A_{4}$ \vspace{.25cm} \\
		 E-$6:$ $x_1 \sim x_2$, for $x_1 \in A_{2}, x_2=e$ & E-$13:$ $x_1 \sim x_2$, for $x_{1},  x_{2} \in A_{4}, x_2 = r^{2^{k-1}p}$ \vspace{.25cm} \\
	     E-$7:$ $x_1 \sim x_2$, for $x_1, x_2 \in A_{5}$ &  \\
		& \\
		\hline	
	\end{tabular}
\end{center}

Note that,  E-$2$ and E-$10$ defines same type of edges,  E-$3$ and E-$6$ defines same type of edges. So, here we consider E-$6$ and E-$10$. Also, edges of types E-$4$, E-$5$, E-$7$, E-$8$, E-$9$, E-$12$, and E-$13$ are subsets of E-$1$ type edges. Let $M_{j}^{i}$ is the number of matchings of order $i$ of type $M_j$ for $1 \leq j \leq 15$.

\begin{center}
	
	\noindent \textbf{Table 3} : Types of matchings
	
	\renewcommand{\arraystretch}{.9}
	\begin{tabular}{|p{5cm}|p{6.5cm}|p{8cm}|}
		
		\hline
		Type &  Type   \vspace{.25cm} \\  \hline
        $M_{1}:$  Between $E$-$1$  & $M_{9}:$  Between $E$-$6$, $E$-$10$ \vspace{.25cm} \\
	    $M_{2}:$  Between $E$-$6$ & $M_{10}:$  Between $E$-$6$, $E$-$11$ \vspace{.25cm} \\
	    $M_{3}:$  Between $E$-$10$ & $M_{11}:$  Between $E$-$1, $E$-$10$, E$-$11$ \vspace{.25cm} \\
	    $M_{4}:$  Between $E$-$11$ & $M_{12}:$  Between $E$-$6$, $E$-$4$, $E$-$11$ \vspace{.25cm} \\
	    $M_{5}:$  Between $E$-$1$, $E$-$10$ & $M_{13}:$  Between $E$-$6$, $E$-$10$, $E$-$11$ \vspace{.25cm} \\
	    $M_{6}:$  Between $E$-$1$, $E$-$11$ & $M_{14}:$  Between $E$-$6$, $E$-$4$, $E$-$10$ \vspace{.25cm} \\ 
	    $M_{7}:$  Between $E$-$10$, $E$-$11$ & $M_{15}:$  Between $E$-$6$, $E$-$4$, $E$-$10$, $E$-$11$ \vspace{.15cm} \\
	    $M_8:$  Between $E$-$6$, $E$-$4$ &   \\
	    \hline
	\end{tabular}
\end{center}

\smallskip 
The following types of distinct matchings occurs among the edges of $\mathcal{P}(G)$.
\begin{description}

	\item [$M_{1}$ Type:] Since the subgraph induced by $A_{1}$ is $K_{2^kp}$, the number of $M_1$ type matchings is the matchings in $K_{2^kp}$. Table 1 lists the matchings in $K_{2^kp}$, where $m_i$ is the number of matchings of order $i$ for $1 \leq i \leq {2^{k-1}p}.$ Thus, $M^i_1 = \frac{1}{i} \prod_{k=0}^{i-1}\binom{2^kp-2k}{2}$ for  $1 \leq i \leq {2^{k-1}p}.$
	
	\item [$M_{2}$ Type:] There are $2^{k-1}p$ edges of type $E$-$6$ which provide $2^{k-1}p$  matchings of order 1. Thus, $M^{1}_2 = 2^{k-1}p$.
	
	\item [$M_{3}$ Type:] There are $2^{k}p$ edges of type $E$-$10$, which provide $2^{k}p$ matchings of order 1. So, $M_{3}^{1}=2^{k}p.$
	
	Suppose $v_1 \sim v_2$ is an $E$-$10$ type edge with $v_{1} \in A_{6}^{j}$ for a fixed $0 \leq j \leq 2^{k-2}p-1$ and $v_{2} \in \Omega$.
	Then, every edge $u_1 \sim u_2$ of $E$-$10$ type with $u_1 \in A_{6} \setminus \{ v_{1}\}$ and $u_2 \in \Omega \setminus \{ v_{2}\}$ together with edge $v_1 \sim v_2$ forms a matching of order $2$.
	So, $$M_{3}^{2} = 2 \cdot \frac{2^{k-1}p(2^{k-1}p-1)}{2}= 2^{k-1}p(2^{k-1}p-1).$$

	\item [$M_{4}$ Type:]  Since there are $2^{k-2}p$ distinct edges of $E$-$11$ type, so for every $i$ we get a matching $M_{4}^{i} = \binom{2^{k-2}p}{i}$ for $1 \leq i \leq 2^{k-2}p$.
	
	\item [$M_{5}$ Type:] This type of matching is obtained by taking matchings of order $1$ and $2$ between the edges of type $E$-$10$ and a matching of order $\ell$ between the edges of $E$-$1$ type having end vertices other than $e$ or $r^{2^{k-1}p}$, where $1 \leq \ell \leq 2^{k-1}p$.
	
	 Now, there are two possibilities for order $i$ matching as follows:
	
	$(i)$ Suppose $v_1 \sim v_2$ is an $E$-$10$ type edge with $v_{1} \in A_{6}^{j}$ for a fixed $0 \leq j \leq 2^{k-2}p-1$ and $v_{2} \in \Omega$.
	Then, every edge $u_1 \sim u_2$ of $E$-$1$ type except $u_1 \sim v_2$ together with edge $v_1 \sim v_2$ forms a matching. Thus, in this case, every matching of order $i$ is obtained by adding one edge of $E$-$10$ type to a matching of order $i-1$ in $K_{2^kp-1}$. By Table 1, there are $\frac{1}{i-1} \prod \limits_{s=0}^{i-2} \binom{2^kp-2s-1}{2}$ matchings of order $i-1$ in $K_{2^kp-1}.$  Since there are $2^{k}p$ distinct edges of $E$-$10$ type, for $2 \leq i \leq 2^{k-1}p$, we get
	
	$$M_{5}^{i} = 2^{k}p \cdot \frac{1}{i-1} \prod_{s=0}^{i-2} \binom{2^kp-2s-1}{2}.$$
	
	$(ii)$ Consider the edges $v_1 \sim e$ and $v_2 \sim r^{2^{k-1}p}$ of $E$-$10$ type with $v_{1}, v_{2} \in A_{6}^{j}$ for $0 \leq j \leq 2^{k-2}p-1$. Then the edges $v_1 \sim e$ and $v_2 \sim r^{2^{k-1}p}$ together with edges $u_1 \sim u_2$ of $E$-$1$ type such that $u_1, u_2 \notin \Omega$ forms a matching. Thus, in this case every matching of order $i$ is obtained by adding a matching of order 2 between $E$-$10$ type edges to a matching of order $i-2$ in $K_{2^kp-2}$. Note that, there are $\frac{1}{i-2} \prod\limits_{s=0}^{i-3} \binom{2^kp-2s-2}{2}$ matchings of order $i-2$ in $K_{2^kp-2}.$ Also, there are $\frac{2^{k-1}p(2^{k-1}p-1)}{2}$ matchings of order 2 between $E$-$10$ type edges.  Thus, for $2 \leq i \leq 2^{k-1}p+1$, we get
	
	$$M_{5}^{i} = \frac{2^{k-1}p(2^{k-1}p-1)}{2}  \cdot \frac{1}{i-2} \prod_{s=0}^{i-3} \binom{2^kp-2s-2}{2}.$$

	Thus, for $2 \leq i \leq 2^{k-1}p$, $$M_{5}^{i} = 2^{k}p \cdot \frac{1}{i-1} \prod_{s=0}^{i-2} \binom{2^kp-2s-1}{2} + \frac{2^{k-1}p(2^{k-1}p-1)}{2}  \cdot \frac{1}{i-2} \prod_{s=0}^{i-3} \binom{2^kp-2s-2}{2}$$
	and 
	$$M_{5}^{2^{k-1}p+1} = \frac{2^{k-1}p(2^{k-1}p-1)}{2}  \cdot \frac{1}{2^{k-1}p-1} \prod_{s=0}^{2^{k-1}p-2} \binom{2^kp-2s-2}{2}.$$
	
	\item [$M_{6}$ Type:] As edges of type $E$-$1$ and type $E$-$11$ are disjoint, so a matching of type $M_6$ is obtained by including a matching of $E$-$1$ type edges to a matching of $E$-$11$ type edges. Note that, every edges of $E$-$1$ type is also an edge of $K_{2^kp}$ and this gives $F_{6}^{\ell}$ matchings of order $\ell$ between them, where $1 \leq \ell \leq 2^{k-1}p.$ Since there are $2^{k-2}p$ distinct edges of $E$-$11$ type, so we get $H^m_6$ matchings of order $m$ between them, where $1 \leq m \leq 2^{k-2}p$. This gives $H_{6}^{m} = \binom{2^{k-2}p}{m}.$ Thus, for $2 \leq i \leq 3.2^{k-2}p,$ we have
	
	\[\begin{aligned}
	M_{6}^{2} &= F_{6}^{1}H_{6}^{1}, \\
	M_{6}^{3} &= F_{6}^{1}H_{6}^{2} + F_{6}^{2}H_{6}^{1}, \\
	M_{6}^{4} &= F_{6}^{1}H_{6}^{3} + F_{6}^{2}H_{6}^{2} + F_{6}^{3}H_{6}^{1}, \\
	& ~~ \vdots \\
	M_{6}^{i} &= \sum\limits_{\jmath=1}^{i-1} F_{6}^{\jmath}H_{6}^{i-\jmath}, \text{ where } F_{6}^{\jmath} = 0 \text{ for } j > 2^{k-1}p \text{ and } H_{6}^{i-\jmath} = 0 \text{ for } i-\jmath > 2^{k-2}p.
	\end{aligned}\]

	\item[$M_{7}$ Type:] Every matching of type $M_7$ is obtained by taking matchings of order $1$ and $2$ between the edges of $E$-$10$ type and a matching of order $\ell$ between the edges of $E$-$11$ type, where $1 \leq \ell \leq 2^{k-2}p-1$.  This gives,
	\[M_{7}^{2}= 2^kp \binom{2^{k-2}p-1}{1}= 2^kp(2^{k-2}p-1),\]
	
	and for $3 \leq i \leq 2^{k-2}p-1,$ we have
	
	\[M_{7}^{i}=2^{k}p \cdot \binom{2^{k-2}p-1}{i-1} + 2.2^{k-2}p \cdot \binom{2^{k-2}p-1}{i-2} + 2 \frac{(2^{k-1}p)(2^{k-1}p-2)}{2} \binom{2^{k-2}p-2}{i-2}.\]
	
	\item[$M_{8}$ Type:] 
	A matching of type $M_8$ is obtained by taking one edge of $E$-$6$ type and at least one edge of $E$-$4$ type.
	Note that, every edges of $E$-$4$ type is also an edge of $K_{2^kp-1}$, so every matching among edges of $E$-$4$ is also a matching of $K_{2^kp-1}$. Hence, a matching of type $M_8$ with order $i$ is obtained by adding one edge of type $E$-$6$ to a matching of order $i-1$ in $K_{2^kp-1}$.
	By Table 1, there are $\frac{1}{i-1} \prod \limits_{s=0}^{i-2} \binom{2^kp-2s-1}{2}$ matchings of order $i-1$ in $K_{2^kp-1}$. Since there are $2^{k-1}p$ edges of $E$-$6$ type, the number of matchings $M^i_8$ of order $i$ for $2 \leq i \leq 2^{k-1}p$ is
	$$\frac{2^{k-1}p}{i-1} \prod_{s=0}^{i-2} \binom{2^kp-2s-1}{2}.$$
	
    \item[$M_{9}$ Type:]  Suppose $v_1 \sim v_2$ is an $E$-$10$ type edge with $v_{1} \in B_{6}^{j}$, for a fixed $0 \leq j \leq 2^{k-2}p-1$ and $v_{2} \in \Omega$.
    Then, every edge $u_1 \sim u_2$ of type $E$-$10$ with $u_1 \in B_{6}$ and $u_2 \in \Omega \setminus \{ e\}$ together with a matching of order 1 between the edges of $E$-$6$ type gives
     a matching of order $2$.
    So, $$M_{9}^{2} = 2^{k-1}p \cdot 2^{k-1}p = 2^{2k-2}p^2.$$
	
	\item[$M_{10}$ Type:] Since edges of types $E$-$6$ and $E$-$11$ are disjoint, so such a matching is obtained by including a matching among $E$-$6$ type edges to a matching of $E$-$11$ type edges. So, for $2 \leq i \leq 2^{k-2}p+1,$	$$M_{10}^{i}= 2^{k-1}p \cdot \binom{2^{k-2}p}{i-1}.$$
	
	\item[$M_{11}$ Type:] A matching of such type is obtained by taking at least one edge of $E$-$1$ type, at least one edge of $E$-$10$ type, and at least one edge of $E$-$11$ type. 
	Then, we have the following cases:
	
	\noindent{\textbf{Case 1 :}} Consider an order one matching between edges of type $E$-$10$ and let $B_{11}^{1}, C_{11}^{\ell}$, and $D_{11}^{m}$ denote the number of matchings of order $1, \ell,$ and $m$ among the edges of types $E$-$10$, $E$-$1$, and $E$-$11$ respectively. Then an $i$ order matching is obtained by taking a $\ell$ order matching among $E$-$1$ type edges having end vertices other than $e$ or $r^{2^{k-1}p}$ and $i-\ell-1$ order matching among $E$-$11$ type edges. Note that in this case $1 \leq \ell \leq 2^{k-1}p-1$, $1 \leq m \leq 2^{k-2}p-1$, and $B_{11}^{1} = 2^kp, \ C^{\ell}_{11} = \frac{1}{\ell} \prod_{s=0}^{\ell-1} \binom{2^kp-2s-1}{2}$, $D_{11}^{m} = \binom{2^{k-2}p-1}{m}$. Then, for $ 3 \leq i \leq 3 \cdot 2^{k-2}p-1$, we have
	
	\[\begin{aligned}
	N_{11}^{3} &= B_{11}^{1}C_{11}^{1}D_{11}^{1}, \\
	N_{11}^{4} &= B_{11}^{1}C_{11}^{1}D_{11}^{2} + B_{11}^{1}C_{11}^{2}D_{11}^{1}, \\ & ~~ \vdots \\
	N_{11}^{i} &= \sum\limits_{\jmath=1}^{i-2} B_{11}^{1}C_{11}^{\jmath}D_{11}^{i-\jmath-1}.
	\end{aligned}\]
	
	$\text{ where } C_{11}^{\jmath} = 0 \text{ for } j > 2^{k-1}p -1 \text{ and } D_{11}^{i-\jmath-2} = 0 \text{ for }  i-\jmath-1 > 2^{k-2}p-1.$
	
	\noindent{\textbf{Case 2 :}} If we take a matching of order two among the edges of type $E$-$10$, then we have the following sub cases:
	
	$(i)$ Suppose one of the end vertices of both the edges of the matching among edges of type $E$-$10$ belong to $A^j_6$ for a fix $j$ and let $R_{11}^{2}, S_{11}^{\ell}$, and $T_{11}^{m}$ denote the number of matchings of order $2, \ell,$ and $m$  among the edges of types $E$-$10$, $E$-$1$, and $E$-$11$ respectively. Then, in this case, $1 \leq \ell \leq 2^{k-1}p-1$, $1 \leq m \leq 2^{k-2}p-1$, and $R_{11}^{2} = 2^{k-1}p, \ S^{\ell}_{11} = \frac{1}{\ell} \prod_{s=0}^{\ell-1} \binom{2^kp-2s-2}{2}$, $T_{11}^{m} = \binom{2^{k-2}p-1}{m}$. Thus, for $ 4 \leq i \leq 3 \cdot 2^{k-2}p$, we have
	\[\begin{aligned}
	P_{11}^{4} &= R_{11}^{2}S_{11}^{1}T_{11}^{1}, \\
	P_{11}^{5} &= R_{11}^{2}S_{11}^{1}T_{11}^{2} + R_{11}^{2}S_{11}^{2}T_{11}^{1}, \\
	P_{11}^{6} &= R_{11}^{2}S_{11}^{1}T_{11}^{3} + R_{11}^{2}S_{11}^{2}T_{11}^{2} + R_{11}^{2}S_{11}^{3}T_{11}^{1}, \\
	& ~~ \vdots \\
	P_{11}^{i} &= \sum\limits_{\jmath=1}^{i-3} R_{11}^{2}S_{11}^{\jmath}T_{11}^{i-\jmath-2}.
	\end{aligned}\]
	
	$\text{ where } S_{11}^{\jmath} = 0 \text{ for } j > 2^{k-1}p -1 \text{ and } T_{11}^{i-\jmath-2} = 0 \text{ for }  i-\jmath -2 > 2^{k-2}p-1.$
	
	$(ii)$ Suppose one of the end vertices of both the edges of the matching among edges of type $E$-$10$ belong to different $A^j_6$ and $X_{11}^{2}, Y_{11}^{\ell}$, and $Z_{11}^{m}$ denote the number of matchings of order $2, \ell,$ and $m$  among the edges of types $E$-$10$, $E$-$1$, and $E$-$11$ respectively. Then, in this case $1 \leq \ell \leq 2^{k-1}p-1$, $1 \leq m \leq 2^{k-2}p-2$, and $X_{11}^{1} = 2 \cdot \frac{2^{k-1}p(2^{k-1}p-1)}{2}, \ Y^{\ell}_{11} = \frac{1}{\ell} \prod_{s=0}^{\ell-1} \binom{2^kp-2s-2}{2}$, $Z_{11}^{m} = \binom{2^{k-2}p-2}{m}$. Thus, for $ 4 \leq i \leq 3 \cdot 2^{k-2}p -1$, we have
	
	\[\begin{aligned}
	Q_{11}^{4} &= X_{11}^{2}Y_{11}^{1}Z_{11}^{1}, \\
	Q_{11}^{5} &= X_{11}^{2}Y_{11}^{1}Z_{11}^{2} + X_{11}^{2}Y_{11}^{2}Z_{11}^{1}, \\
	Q_{11}^{6} &= X_{11}^{2}Y_{11}^{1}Z_{11}^{3} + X_{11}^{2}Y_{11}^{2}Z_{11}^{2} + X_{11}^{2}Y_{11}^{3}Z_{11}^{1}, \\
	& ~~ \vdots \\
	Q_{11}^{i} &= \sum\limits_{\jmath=1}^{i-3} X_{11}^{2}Y_{11}^{\jmath}Z_{11}^{i-\jmath-2}.
	\end{aligned}\]
	
	$\text{ where } Y_{11}^{\jmath} = 0 \text{ for } j > 2^{k-1}p -1 \text{ and } Z_{11}^{i-\jmath-2} = 0 \text{ for }  i-\jmath -2 > 2^{k-2}p-2.$

	Thus, $M^{3}_{11} = N^{3}_{11}$, for $4 \leq i \leq 3 \cdot 2^{k-2}p -1$, $M^{i}_{11} = N^i_{11} + P^{i}_{11} + Q^{i}_{11}$ and $M^{3 \cdot 2^{k-2}p}_{11} = Q^{3 \cdot 2^{k-2}p}_{11}.$
	
	\item[$M_{12}$ Type:]  A matching of this type is obtained by taking one edge of type $E$-$6$, at least one edge of type $E$-$4$, and at least one edge of type $E$-$11$. Note that in this case, edges of types $E$-$6$, $E$-$4$, and $E$-$11$ are disjoint and every edges of type $E$-$4$ is also an edge of $K_{2^kp-1}$. Let $I_{12}^{1}$ denote the matchings of order $1$ between the edges of type $E$-$6$, $J_{12}^{\ell}$ denote the matchings of order $\ell$ between the edges of type $E$-$4$, where $1 \leq \ell \leq 2^{k-1}p-1$, and $L_{12}^{m}$ denote matchings of order $m$ between the edges of type $E$-$11$, where $1 \leq m \leq 2^{k-2}p$. Clearly, $I^{1}_{12} = 2^{k-1}p, \ J^{\ell}_{12} = \frac{1}{\ell} \prod \limits_{s=0}^{\ell-1} \binom{2^kp-2s-1}{2}, L^{m}_{12} = \binom{2^{k-2}p}{m}.$ So, for $3 \leq i \leq 3 \cdot 2^{k-2}p$, the number of matchings $M^i_{12}$ of order $i$ is
	
	\[\begin{aligned}
	M_{12}^{3} &= I_{12}^{1}J_{12}^{1}L_{12}^{1}, \\
	M_{12}^{4} &= I_{12}^{1}J_{12}^{1}L_{12}^{2} + I_{12}^{1}J_{12}^{2}L_{12}^{1}, \\
	M_{12}^{i} &= \sum\limits_{\jmath=1}^{i-2} I_{12}^{1}J_{12}^{\jmath}L_{12}^{i-\jmath-1},
	\end{aligned}\]
	
	$\text{ where } J_{12}^{\jmath} = 0 \text{ for } j > 2^{k-1}p -1 \text{ and } L_{12}^{i-\jmath-1} = 0 \text{ for }  i-\jmath -2 > 2^{k-2}p.$
	
%
%

	\item[$M_{13}$ Type:] 
	A matching of such type is obtained by taking one edge of type $E$-$6$ together with one edge of type $E$-$10$ and at least one edge of type $E$-$11$. Since the vertex $e$ is an end vertex of $E$-$6$ type edges, so the vertex $r^{2^{k-1}p}$ is an end vertex of $E$-$10$ type edge.  Let $O_{13}^{1}$ denote the matchings of order $1$ between the edges of type $E$-$6$, $P_{13}^{1}$ denote  matchings of order $1$ between the edges of type $E$-$10$ having end vertices other than $e$, and $Q_{13}^{m}$ denote matchings of order $m$ between edges of type $E$-$11$, where $1 \leq m \leq 2^{k-2}p-1$. Clearly, $O^{1}_{13} = 2^{k-1}p, \ P^{1}_{13} = 2^{k-1}p, Q^{m}_{13} = \binom{2^{k-2}p-1}{m}.$ So, for $3 \leq i \leq 2^{k-2}p + 1$, the number of matchings $M^i_{13}$ of order $i$ is

	\[\begin{aligned}
	M_{13}^{i} & = O_{13}^{1}P_{13}^{1}Q_{13}^{i-2}.
	\end{aligned}\]

	\item[$M_{14}$ Type:] A matching of this type is obtained by taking an edge of type $E$-$6$ together with an edge of type $E$-$10$ and at least one edge of type $E$-$4$. 
	Since the vertex $e$ is an end vertex of $E$-$6$ type edges, so the vertex $r^{2^{k-1}p}$ is an end vertex of $E$-$10$ type edges. This gives that a matching between the edges of type $E$-$4$ having end vertices other than $r^{2^{k-1}p}$ is a matching of $K_{2^kp-2}$. Let $L_{14}^{1}$ denote a matching of order $1$ between the edges of $E$-$6$ type,  $M_{14}^{1}$ denote a matching of order $1$ between the edges of $E$-$10$ type, and $N_{14}^{1}$ denote a matching of order $\ell$ in $K_{2^kp-2}$, where $1 \leq \ell \leq 2^{k-1}p-1$. Clearly, $L^{1}_{14} = 2^{k-1}p, \ M^{1}_{14} = 2^{k}p, \ N^{\ell}_{14} = \frac{1}{\ell} \prod \limits_{s=0}^{\ell-1} \binom{2^kp-2s-2}{2}.$ So, for $3 \leq i \leq 2^{k-1}p+1$, the number of matchings $M^i_{14}$ of order $i$ is
	
	\[\begin{aligned}
	M_{14}^{i} & = L_{14}^{1}M_{14}^{1}N_{14}^{i-2}.
	\end{aligned}\]

	\item[$M_{15}$ Type:] A matching of $M_{15}$ type is obtained by taking an edge of type $E$-$6$ together with one edge of type $E$-$10$, at least one edge of $E$-$4$ type, and at least one edge of $E$-$11$ type. 
	Since $e$ is an end vertex of $E$-$6$ type edges, so the vertex $r^{2^{k-1}p}$ is an end vertex of $E$-$10$ type edges. This gives that every matching between the edges of type $E$-$4$ having end vertices other than $r^{2^{k-1}p}$ is a matching of $K_{2^kp-2}$. Let $R_{15}^{1}$ denote a matching of order $1$ between the edges of type $E$-$6$,  $S_{15}^{1}$ denote a matching of order $1$ between the edges of $E$-$10$ type, $T_{15}^{1}$ denote a matching of order $\ell$ in $K_{2^kp-2}$, where $1 \leq \ell \leq 2^{k-1}p-1$, and $U_{15}^{m}$ denote a matching of order $m$ between edges of type $E$-$11$, where $1 \leq m \leq 2^{k-2}p-1$. Clearly, $R^{1}_{15} = 2^{k-1}p, \ S^{1}_{15} = 2^{k}p, \ T^{\ell}_{15} = \frac{1}{\ell} \prod \limits_{s=0}^{\ell-1} \binom{2^kp-2s-2}{2}, U^{m}_{15} = \binom{2^kp-1}{m}.$ So, for $4 \leq i \leq 3 \cdot 2^{k-2}p$, the number of matchings $M^i_{15}$ of type order $i$ is

	\[\begin{aligned}
	M_{15}^{4} &= R_{15}^{1}S_{15}^{1}T_{15}^{1}U_{15}^{1}, \\
	M_{15}^{5} &= R_{15}^{1}S_{15}^{1}T_{15}^{1}U_{15}^{2} + R_{15}^{1}S_{15}^{1}T_{15}^{2}U_{15}^{1}, \\
	M_{15}^{6} &= R_{15}^{1}S_{15}^{1}T_{15}^{2}U_{15}^{2} + R_{15}^{1}S_{15}^{1}T_{15}^{3}U_{15}^{1}, \\
	& ~~ \vdots \\
	M_{15}^{i} &= \sum\limits_{\jmath=1}^{i-3} R_{15}^{1}S_{15}^{1}T_{15}^{\jmath}U_{15}^{i-\jmath-2},
	\end{aligned}\]

 $\text{ where } T_{15}^{\jmath} = 0 \text{ for } j > 2^{k-1}p -1 \text{ and } U_{15}^{i-\jmath-2} = 0 \text{ for }  i-\jmath -2 > 2^{k-2}p-1.$
\end{description} This proves the theorem. \hfill $\Box$



\section*{Acknowledgment}

The first author is thankful to the Ministry of Education, New Delhi, India for financial support.

\end{document}